\documentclass[onefignum,onetabnum]{siamart190516}

\usepackage{tikz}
\usetikzlibrary{scopes}
\usetikzlibrary{arrows}
\usetikzlibrary{arrows.meta}
\usetikzlibrary{decorations.pathreplacing,decorations.markings}
\usepackage{lipsum}
\usepackage{amsfonts}
\usepackage{amsmath}
\usepackage{graphicx}
\usepackage{epstopdf}
\usepackage{algorithmic}
\usepackage{centernot}
\usepackage{color}
\ifpdf
\DeclareGraphicsExtensions{.eps,.pdf,.png,.jpg}
\else
\DeclareGraphicsExtensions{.eps}
\fi

\newsiamremark{remark}{Remark}
\newsiamremark{hypothesis}{Hypothesis}
\crefname{hypothesis}{Hypothesis}{Hypotheses}
\newsiamthm{claim}{Claim}

\headers{Numerical procedure for hybrid systems}{R. Pytlak and D. Suski}

\title{Numerical procedure for optimal control of hybrid systems with sliding modes, Part I}
\author{Rados\L aw Pytlak\thanks{Faculty of Mathematics and Information Science,
		Warsaw University of Technology, 00-665 Warsaw, Poland 
		(\email{r.pytlak@mini.pw.edu.pl}),}
	\and Damian Suski\thanks{Institute of Automatic Control and Robotics,
		Warsaw University of Technology, 02-525 Warsaw, Poland 
		(\email{d.suski@mchtr.pw.edu.pl}).}
}

\usepackage{amsopn}

\begin{document}

\maketitle

\begin{abstract}
This paper concerns the numerical procedure for solving hybrid optimal control problems with sliding modes. The proposed procedure has several features which distinguishes it from the other procedures for the problem. First of all a sliding mode is coped with differential-algebraic equations (DAEs) and that guarantees accurate tracking of the sliding motion surface. The second important feature is the calculation of cost and constraints functions gradients with the help of adjoint equations. The adjoint equations presented in the paper take into account sliding motion and exhibit jump conditions at transition instants. 
The procedure uses the discretization of system equations by Radau IIA Runge--Kutta scheme and the evaluation of optimization functions gradients with the help of the adjoint equations stated for  discretized system equations. 
In the first part of the paper we demonstrate the correspondence between the discrete adjoint equations and the discretized version of the continuous adjoint equations in the case of system equations described by ODEs. We show that the discrete adjoint state trajectories converge to their continuous counterparts in the case of ODEs. 
\end{abstract}

% REQUIRED

% REQUIRED

\section{Introduction}

The systems with sliding modes considered in our paper state a special class of so called hybrid systems (\cite{ss2000}). Many technical and non-technical systems exhibit the hybrid discrete-continuous behavior, therefore optimal control problems for hybrid systems were considered in many papers (\cite{w1966}, \cite{bh1975},  \cite{s1999}, \cite{s2004}, \cite{sc2007}, \cite{sc2009}, \cite{tc2010}, \cite{tc2011}, \cite{tc2013}, \cite{pc2013}). Some of them discussed also algorithms for these problems (\cite{s2004}, \cite{sc2007}, \cite{sc2009}, \cite{tc2010}, \cite{tc2013}). These algorithms refer to the strong version of the maximum principle for hybrid systems. It means that their applicability is limited since the class of problems which can be solved by algorithms using locally the strong maximum principle is much smaller than the class of control problems which could be tackled by methods which invoke locally the weak maximum principle condition.

The algorithm based on weak variations and aimed at solving optimal control problems with hybrid systems was stated in  (\cite{ps17}). The numerical procedure discussed therein calls the implicit Runge--Kutta method to evaluate system equations and calculates gradients of the optimization problem functionals with the help of adjoint equations. The adjoint equations are figured out on the basis of discrete time system equations which are the result of the numerical integration of system equations. It means that the discretization of the adjoint equations is consistent with the discretization of the system equations. The paper presents the results of employing the method to solving some optimal control problems described by hybrid systems expressed by higher index DAEs. However, the paper does not present the convergence analysis of the procedure.

The first step towards the missing convergence analysis is done in the accompanying paper (\cite{ps19}). The optimal control problem discussed in the paper includes, as the special case, the problem analyzed in (\cite{ps17}) provided that the hybrid system evaluates according to ordinary differential equations (ODEs). On the other hand, the hybrid system examined in (\cite{ps19}) is more general since it includes the sliding motion. The paper (\cite{ps19}) states necessary optimality conditions in the form of the weak maximum principle for optimal control problems with hybrid systems treated in such a general setting. It also proposes a globally convergent algorithm for the considered optimal control problems.

The aim of this paper is to show that a conceptual algorithm stated in (\cite{ps19}) can be a basis for building its efficient implementable version. It can be achieved by following the approach stated in (\cite{ps17}). Since the system movement in a sliding mode is represented by DAEs with index 2 the numerical procedure from (\cite{ps17}) is applicable to these hybrid systems. Therefore, the gradients of the functions defining the optimal control problem can be evaluated using the discrete adjoint equations consistent with the discretized system equations. However, in contrast to the results of (\cite{ps17}) we show that the gradients  elicited in that way converge to the gradients of functionals defining the continuous time version of the optimal control problem provided that the step sizes in the integration procedure approach zero.  

The computational approach to hybrid systems we advocate in this paper (and in \cite{ps17},\cite{ps19},\cite{ps19b}) assumes that for a given control function we attempt to follow true systems trajectories as close as possible by using integration procedures with a high order of convergence. It is in contrast to the approach which uses approximations to the discontinuous right-hand side for the differential equations or inclusions by a smooth right--hand side (\cite{sa2010},\cite{sa2012},\cite{mhrl2015}). The smoothing approach demonstrates that it can be used to find approximate solutions to optimal control problems with systems which exhibit sliding modes however there is still a lack of evidence that it is capable of efficiently finding optimal controls and state trajectories with high accuracy by choosing proper values of smoothing parameters.

The organization of our paper is as follows. In the next section we describe the systems with sliding motion. Then, we state continuous adjoint equations for these systems. In Section \ref{SecDiscrete} we formulate the discretized system equations and the discrete adjoint equations and show the correspondence between the discrete adjoint equations and discretization of the continuous adjoint equations. On that basis, we show that the discrete adjoint trajectories converge to their continuous counterparts for step sizes approaching zero. 
In the first part of the paper we give an estimate for the order of convergence of discrete adjoint equations assuming that a hybrid system is described by ODEs. In the second part of the paper (\cite{ps19b}) we carry out analogical analysis for the case when system equations correspond to index 2 DAEs. The second part of the paper complements the result of the first part by paying attention to the jumps of adjoint variables which occur when the hybrid system changes its discrete state. We show that the jumps in the discrete adjoint equations converge to jumps of their continuous counterparts provided that the integration step sizes approach zero. In the second part of the paper we also present the numerical results of applying our numerical procedure to several optimal control problems.  

\section{Sliding mode}

The paper concerns a numerical procedure for optimal control problems described by hybrid systems with sliding modes. In this section we present the concept of the sliding motion. Let us consider a hybrid system, which evolves according to ODEs
\begin{equation}
x' = f^1(x,u), \label{SlidingMEq1}
\end{equation}
if the condition
\begin{equation}
g(x)< 0 \label{SlidingMInv1}
\end{equation}
is satisfied. The system evolves according to ODEs
\begin{eqnarray}
x' = f^2(x,u) \label{SlidingMEq2}
\end{eqnarray}
if the opposite condition
\begin{equation}
g(x)> 0 \label{SlidingMInv2}
\end{equation}
holds.
Suppose that the system first evolves according to (\ref{SlidingMEq1}) and at the moment $ t_t $ reaches the switching surface
\begin{equation}
\Sigma = \{ x\in \mathbb{R}^n: g(x) = 0 \}. 
\end{equation}
If at $t_t$ the following conditions are satisfied 
\begin{eqnarray}
g_x(x(t_t))f^1(x(t_t),u(t_t)) & > & 0 \label{SlidingMSw12a} \\
g_x(x(t_t))f^2(x(t_t),u(t_t)) & >  & 0. \label{SlidingMSw12b}
\end{eqnarray}
the system starts to evolve according to ODEs (\ref{SlidingMEq2}). However if at $ t_t $ the condition 
\begin{eqnarray}
g_x(x(t_t))f^2(x(t_t),u(t_t)) & < & 0 \label{SlidingMSwS}
\end{eqnarray}
is satisfied instead of (\ref{SlidingMSw12b}), then both vector fields $ f^1(x,u) $ and $ f^2(x,u) $ point towards the surface $ \Sigma $ and we face the {\it sliding motion} phenomenon (\cite{dl2009}). The system state remains on the surface $\Sigma$ and we say that the system is in the sliding mode, which can be considered as the separate discrete state of the system.

The system motion at the sliding mode can be handled by 
by applying the concept of {\it Filippov solutions} (\cite{dl2009}). We say that a continuous state trajectory is a Filippov solution of the considered hybrid system if
\begin{equation}
x' = \left\{
\begin{array}{ll}
f^1(x,u) & \text{if } g(x) < 0 \\
f_F(x,u) & \text{if } g(x) = 0 \\
f^2(x,u) & \text{if } g(x) > 0 
\end{array}
\right.
\end{equation}
where $ f_F(x,u)$ is a convex combination of $ f^1(x,u) $ and $ f^2(x,u) $
\begin{equation}
f_F(x,u) = (1-\alpha)f^1(x,u) + \alpha f^2(x,u),\ \alpha \in [0,1]. \label{fFDef}
\end{equation} 
Since the coefficient $\alpha$ should guarantee that the system stays on the surface $\Sigma$, the condition
\begin{equation}
g_x(x)f_F(x,u)=0.\label{ffTang}
\end{equation} 
must be satisfied. Therefore, $\alpha$ is a function of $x$ and $u$ and can be evaluated as
\begin{equation}
\alpha(x,u) = \frac{g_x(x)f^1(x,u)}{g_x(x)\left(f^1(x,u)-f^2(x,u) \right)}.\label{AlphaF}
\end{equation} 

The hybrid system leaves the sliding mode when one of the vector fields $ f^1(x,u) $ and $ f^2(x,u) $ points outside the surface $ \Sigma $ at some time $t_t$. This implies that a time $ t_t $ we either have 
\begin{eqnarray}
\alpha(x(t_t),u(t_t)) &=& 0 \label{sw3} \\
g_x(x(t_t))f^1(x(t_t),u(t_t)) & < & 0 \label{sw3a} \\
g_x(x(t_t))f^2(x(t_t),u(t_t)) & < & 0 \label{sw3b}
%\frac{d}{dt} \alpha(x(t_t),u(t_t)) &<&0 \label{sw4}
\end{eqnarray}
and the system  evolves further according to (\ref{SlidingMEq1})  or
\begin{eqnarray}
\alpha(x(t_t),u(t_t)) - 1 &=& 0 \label{sw5} \\
g_x(x(t_t))f^1(x(t_t),u(t_t)) & > & 0 \label{sw5a} \\
g_x(x(t_t))f^2(x(t_t),u(t_t)) & > & 0 \label{sw5b}
%\frac{d}{dt} \alpha(x(t_t),u(t_t)) &>&0 \label{sw6}
\end{eqnarray}
and the system evolves further according to (\ref{SlidingMEq2}).

At the sliding mode the state trajectory should stay in a set $ \Sigma $. During numerical integration we will fulfill that if we integrate the differential-algebraic equations (DAEs) (see \cite{ap1998} for the justification of the use of these equations)
\begin{eqnarray}
x' &=& f_F(x,u) + g_x(x)^Tz \label{slidingDAEDiff} \\
0 &=& g(x) \label{slidingDAEAlg}
\end{eqnarray}
instead of ODEs
\begin{equation}
x' = f_F(x,u).
\end{equation}
One can show that under the condition $g_x(x)\neq 0$ at each $x\in \Sigma$ the equations (\ref{slidingDAEDiff})-(\ref{slidingDAEAlg}) are index two DAEs (\cite{ps19}). 

\section{Optimal control problem with sliding modes } 

The simple (due to the minimal number of discrete states) optimal control problem with a hybrid system exhibiting a sliding mode---{\bf (P)} can be stated as follows:
\begin{eqnarray}
&\min_{u} \phi(x(t_f))& \label{defOptControlProblemCost}
\end{eqnarray}
subject to the constraints
\begin{eqnarray}
\begin{array}{ll}
x' = f_1(x,u) &\ \ \ {\rm if\ }\ g(x)<0  \\
x' = f_2(x,u) &\ \ \ {\rm if\ }\ g(x)>0  \\
\left \{
\begin{array}{lll}
x' & = & f_F(x,u) + g_x(x)^Tz \nonumber \\
0 &=& g(x) \nonumber
\end{array}\right .  &\ \ \ {\rm if\ }\ {\rm at\ a\ sliding\ mode}
%g(x)=0 
\end{array} \label{SystemEq}
\end{eqnarray}
and the terminal constraints
\begin{eqnarray}
&g^1_i(x(t_{f})) = 0\ \forall i\in E & \\
&g^2_j(x(t_{f})) \leq 0\ \forall j\in I.& \label{defOptControlProblemEndIneq}
\end{eqnarray}
Furthermore, due to the fact that $ \alpha(x,u) $ depends on control functions $u$ the admissible control is a function $u:[t_0,t_f] \rightarrow {\cal U}$, where 
\begin{eqnarray}
&{\displaystyle
	{\cal U} =\left \{ u\in {\cal L}^2_m[T]:\ u(t)=\sum_{j=1}^N {\bf 1}_j(t)u^N(j),\ u^N(j)\in
	U\subset  \mathbb{R}^m\right \} }
%\label{l2d}
\nonumber
\end{eqnarray}
Here, 
\begin{eqnarray}
&{\displaystyle {\bf 1}_j(t)  =   \left \{ \begin{array}{ll}
1, & {\rm if}\ t\in[t_{j-1},t_j]\nonumber \\
0, & {\rm if}\ t\not\in [t_{j-1},t_j]\nonumber 
\end{array} \right . , }
\nonumber
\end{eqnarray}
$t_j=jt_f/N$, $j\in \{0,1,\ldots,N\}$, and $U$ is a closed convex set. The definition of the set of admissible control implies that we are dealing with piecewise constant control functions.

Since our control functions are not continuous it is useful to introduce the notation:
\begin{eqnarray}
&{\displaystyle \lim_{t\rightarrow t_t,t<t_t)} u(t) = u(t_t^-), \ \ 
\lim_{t\rightarrow t_t,t>t_t)} u(t) = u(t_t^+),} \nonumber 
\end{eqnarray}
which will be also applied to other functions, for example we will write $x(t_t^-)$, $\lambda(t_t^+)$.
The convention is also applied to  integrals, for example
\begin{eqnarray}
&{\displaystyle \int_{t_0}^{t_t^-}p(x(t),u(t))dt = \lim_{t\rightarrow t_t,t< t_t}\int_{t_0}^t p(x(\tau),u(\tau))d\tau. }\nonumber
\end{eqnarray}

We can show, under some standard assumptions stated in \cite{ps19}, that for given $x(t_0)$ state trajectory $x$ is uniquely defined by $u$ (so we can write $x^u$ instead of $u$) thus the considered control problem can be expressed as an optimization problem over the set of control
functions with the aid of the functions $\bar{F}_0: {\cal U}\rightarrow {\cal R}$,
$\bar{g}^1_i: {\cal U}\rightarrow {\cal R}$ for $i\in E$,
$\bar{g}^2_j: {\cal U}\rightarrow {\cal R}$ for $j\in I$:
\begin{eqnarray}
\bar{F}_0(u) & = & \phi (x^u(t_f)) \nonumber\\
\bar{g}^1_i (u) & = &  g^1_i (x^u(t_f)) \ \forall i\in E\nonumber\\
\bar{g}^2_j (u) & = &  g^2_j (x^u(t_f))\ \forall j\in I.\nonumber
\end{eqnarray}

Eventually, the optimal control problem we want to solve can be stated as 
\begin{eqnarray}
&{\displaystyle
	\min_{u\in {\cal U} } \bar{F}_0 (u)}\label{gen1}
\end{eqnarray}
subject to
\begin{eqnarray}
\bar{g}^1_i(u) & = & 0\ \forall i\in E\label{gen2a}\\
\bar{g}^2_j(u) & \leq & 0\ \forall j\in I.\label{gen2b}
\end{eqnarray}

\section{Calculating the optimal control}
\label{CalcOptCon}

In \cite{ps19} we propose the method for solving the problem ${\bf (P)}$. This is an exact penalty algorithm discussed in \cite{py99} which was used to solve several optimal control problems including those described by DAEs (\cite{py11},\cite{pz14}). Since the optimization method we advocate for solving for solving the problem ${\bf (P)}$  is based on an exact penalty function, 
instead of solving the problem ${\bf (P)}$ we solve the problem ${\bf P_{c}}$
\begin{eqnarray}
&{\displaystyle
	\min_{u\in {\cal U} } \bar{F}_c(u)}
\label{1e}
%\nonumber
\end{eqnarray}
in which the exact penalty function $\{\bar{F}_c(u)$ is defined as follows
\begin{eqnarray}
&{\displaystyle \bar{F}_c(u) = \bar{F}(u) +c \max \left [0,\max_{i\in E} \left |\bar{g}^1_i(u)\right |,\max_{j\in I} \bar{g}^2_j(u)\right ]} 
\label{2e}
\end{eqnarray}

Let $H\in {\cal R}^{mN\times mN}$ be a symmetric matrix such that
\begin{eqnarray}
&{\displaystyle   \nu_1 \| d\|^2 \leq d^THd \leq \nu_2 \| d\|^2 \ \forall d\in {\cal R}^{mN}  }\label{Hcond}
\end{eqnarray}
for some positive $\nu_1$ and $\nu_2$.

For fixed $c$, $H$ satisfying (\ref{Hcond}) and $u$ the direction finding subproblem, ${\bf P_{c}(u)}$, 
for the problem ${\bf (P_c)}$ is:
\begin{eqnarray}
&{\displaystyle
	\min_{d\in {\cal U}-u,\beta\in {\cal R}}\left [\left
	\langle \nabla \bar{F}_0(u),d\right \rangle +c \beta + 1/2 
	d^THd \right ]
} \nonumber  
%\label{Linapp1}
\end{eqnarray}
subject to
\begin{eqnarray}
\left | \bar{g}^1_i(u) +
\left \langle \nabla \bar{g}^1_i(u),d\right \rangle\right |  & \leq &\beta\ \ \forall i\in E
\nonumber \\
\bar{g}^2_j (u) + \left \langle \nabla \bar{g}^2_j(u),d\right \rangle
& \leq & \beta\ \ \forall j\in I.
\nonumber
\end{eqnarray}

The subproblem can be reformulated as an optimization problem over the space ${\cal R}^{mN}$
with the objective function which is strictly convex.
%The problem therefore has the unique solution $(\bar{d},\bar{\beta})$.
%Since this solution depends on $c$, $H$ and $u$, we may 

We define {\it descent function}
$\sigma_{c}^H(u)$ and {\it penalty test function} $t_{c}(u)$ (see \cite{py99} for details),
to be used to test optimality of a control
$u$ and to adjust $c$, respectively, as 
\begin{eqnarray}
&{\displaystyle
	\sigma_{c}^H(u) = \left \langle \nabla \bar{F}_0(u),\bar{d}\right \rangle + 
	c\left [\bar{\beta}-M(u)\right ] }\nonumber
\end{eqnarray}
and                                                          
\begin{eqnarray}
&{\displaystyle
	t_{c}(u) = \sigma_{c}^H(u) + M(u)/c }\label{testf}
\end{eqnarray}
\label{tesfun2}
for given $c> 0$, $H$ and $u\in {\cal U}$. Here, 
\begin{eqnarray}
&{\displaystyle
	M(u) = \max\left [0,\max_{i\in E}\left |\bar{g}^1_i(u)\right |,
	\max_{j\in I}\bar{g}^2_j(u)\right ],}\nonumber
%\label{30h}
\end{eqnarray}

A general algorithm is as follows.
\\[5mm]
{\bf Algorithm}
Fix parameters: $\varepsilon> 0,\ \gamma,\ \eta  \in (0,1)$, $c^0 > 0$, 
$\kappa >1$, matrix $H$.
\begin{enumerate}
	\item Choose the initial control $u_{0}\in {\cal U}$.
	Set $k=0$, $c_{-1}=c^0$.
	\item Let $c_k$ be the smallest number chosen from $\{c_{k-1},\kappa c_{k-1},
	\kappa^2 c_{k-1},\ldots \}$ such that
	the solution $(d_k,\beta_k)$ to the direction finding subproblem  ${\bf P_{c_k}(u_k)}$ 
	satisfies
	\begin{eqnarray}
	&{\displaystyle
		t_{c_k} (u_k) \leq 0. }\label{penadj}
	\end{eqnarray}
	If $\sigma_{c_k}^H(u_k) =0$ then STOP.
	\item Let $\alpha_k$ be the largest number chosen from the set
	$\{1,\eta,\eta^2,\ldots,\}$ such that
	\begin{eqnarray}
	&{\displaystyle
		u_{k+1} = u_k + \alpha_k d_k}\nonumber
	\end{eqnarray}
	satisfies the relation
	\begin{eqnarray}
	\bar{F}_{c_k}(u_{k+1}) - \bar{F}_{c_k}(u_k)
	& \leq &  \gamma\alpha_k\sigma_{c_k}^H ( u_k).\nonumber 
	\end{eqnarray}
	Increase k by one. Go to Step 2.
\end{enumerate}

We can show that the descent function $\sigma_{c_k}^H(u_k)$ is nonpositive valued at each iteration (\cite{py99}).

The main computational effort is associated with the evaluation of function values $\bar{F}_0$, $\bar{g}_i^1$, $i\in E$, $\bar{g}_j^2$, $j\in I$ and their gradients $\nabla \bar{F}_0$, $\nabla\bar{g}_i^1$, $i\in E$, $\nabla\bar{g}_j^2$, $j\in I$. The number of decision variables is equal to $mN$ and for the fine discretization of control variables it can be counted in hundreds. Therefore, it can be more efficient to use adjoint equations instead of relying on sensitivity equations for evaluating gradients (see Chapter 6 in \cite{py99}). 

The cost of evaluating functions values can be significant since it requires the integration of ODEs and DAEs if sliding modes occur. Higher index DAEs usually are integrated by implicit methods (BDF, or Runge--Kutta---\cite{ap1998}, \cite{hlr89}). Since we approximate controls by piecewise constant functions we can expect that multistep integration procedure would have to be often restarted which would deteriorate efficiency of numerical integration. Runge--Kutta methods are one--step methods so they are more appropriate in optimal control applications (\cite{sp98},\cite{py99},\cite{h00},\cite{pp02}). However, we are constrained by implicit methods and that means that the cost of numerical integration is significant due to the necessity of solving nonlinear algebraic equations (and possibly due to evaluating equations functions Jacobians and their factorizations) at each integration step. 

Furthermore, adjoint equations to higher index DAEs are also higher index DAEs, so the mentioned above efficiency issues apply also to evaluating adjoint equations. Since each function gradient requires solution of unique adjoint equations (due to different terminal conditions) at each optimization step we have to solve $|E| + |I| +1$ adjoint equations.

Therefore, it is tempting to use as much as possible information acquired  during system equations integration in numerical evaluation of adjoint equations. It can be achieved, as we will show in the next section, if adjoint equations are built for discrete time system equations resulting from the numerical integration. However, if we do that we have to show whether these discrete time adjoint equations will give the solution to their continuous counterparts provided that integration steps approach zero. This is the main subject of the next section.

The 'reference' continuous time adjoint equations are established in our accompanying paper \cite{ps19}. Let us recall only one case, when  before the switching time $t_t$ the system evolves according to ODEs
\begin{equation}
x' = f^1(x,u),
\label{eqAdjEqsFirstCaseOde}
\end{equation}
at $ t_t $ system enters the sliding mode and then the state trajectory is the solution to DAEs 
\begin{eqnarray}
x' &=& f_F(x,u) + g_x^T(x)z \label{eqAdjEqsFirstCaseDaeDiff} \\
0 &=& g(x) \label{eqAdjEqsFirstCaseDaeAlg}
\end{eqnarray}
up to an ending time $ t_f $. 
At a transition time $ t_t $ the continuous state trajectory meets the switching surface such that the following condition holds
\begin{equation}
g(x(t_t)) = 0. 
\label{eqAdjEqsFirstCaseTrans}
\end{equation} 

The adjoint equations for the system (\ref{eqAdjEqsFirstCaseDaeDiff})--(\ref{eqAdjEqsFirstCaseDaeAlg}) are as follows (\cite{ps19})
\begin{eqnarray}
(\lambda_f^T)'(t) &= &-\lambda_f^T(t) (f_F)_x(x(t),u(t)) \label{eqAdjEqsFirstCaseAdjDaeDiff}  \\
&\ & - \lambda_f^T(t) (g_x^T(x(t))z(t))_x + \lambda_g^T(t) g_x(x(t)) \nonumber \\
0 &= & \lambda_f^T(t) g_x^T(x(t)),\ t\in [t_t^+,t_f]. \label{eqAdjEqsFirstCaseAdjDaeAlg}
\end{eqnarray}
and the consistent terminal values of the adjoint variables $\lambda_f(t_f)$ and $\lambda_g(t_f)$ can be found by solving the following set of equations with respect to the variables $ \lambda_f(t_f), \lambda_g(t_f), \nu_1 $ 
\begin{eqnarray}
\phi_x^T(x(t_f)) + \lambda_f(t_f) & = & \nu_1 g_x^T(x(t_f)) \label{terminal} \\
0 & = & g_x(x(t_f))\lambda_f(t_f)  \\
0 & = & \left(g_x(x(t_f))\right)'\lambda_f(t_f) - \label{slidingDAEIndex2HiddenAdjEndpoint} \\
& & g_x(x(t_f)) \left( f_F \right)_x^T(x(t_f),u(t_f))\lambda_f(t_f) -\nonumber \\
& & g_x(x(t_f)) \left(g_x^T(x(t_f))z(t_f)\right)_x^T\lambda_f(t_f) +  \nonumber \\
& & g_x(x(t_f)) g_x^T(x(t_f)) \lambda_g (t_f) \nonumber
\end{eqnarray}
where $ \nu_1 $ is some real number.

At the transition time $ t_t $ the adjoint variable $ \lambda_f $ undergoes a jump. To calculate the value of $ \lambda_f(t_t^-) $ the following system of equations have to be solved for the variables $ \lambda_f(t_t^-),\ \pi $ 
\begin{eqnarray}
\lambda_f(t_t^-) & = & \lambda_f(t_t^+) - \pi g_x^T(x(t_t))  \\
\lambda_f^T(t_t^-) f_1(x(t_t^-),u(t_t^-)) & = & \lambda_f^T(t_t^+) f_F(x(t_t^+),u(t_t^+)) + \label{odeAdjTransHamilCont} \\
& & \lambda_f^T(t_t^+)g_x^T(x(t_t^+))z(t_t^+) - \lambda_g^T(t_t^+)g(x(t_t^+)). \nonumber
\end{eqnarray}

Eventually we solve adjoint ODEs 
\begin{equation}
(\lambda_f^T)'(t) = - \lambda_f^T(t) (f_1)_x(x(t),u(t)),\ t\in[t_0,t_t^-).
\label{eqAdjEqsFirstCaseAdjOde}
\end{equation}
%on the time interval $[t_0,t_t^-]$.

Now the directional derivative $\langle \nabla \bar{F}_0(u),d\rangle$ of the cost functional $\bar{F}_0$ 
%from (\ref{Linapp1}) 
takes the form 
\begin{eqnarray}
\langle \nabla \bar{F}_0(u),d\rangle &=& -\int_{t_0}^{t_t^-} \lambda_f^T(t) (f_1)_u(x(t),u(t))d(t) dt \label{eqAdjGradient}\\
&&   -\int_{t_t^+}^{t_f} \lambda_f^T(t) (f_F)_u(x(t),u(t))d(t) dt. \nonumber  
\end{eqnarray} 

In the next section we will examine the discrete time versions of equations (\ref{eqAdjEqsFirstCaseAdjOde}) and the approximation to the first term in (\ref{eqAdjGradient}). The approximation to the second term and the analysis of discrete time versions of equations (\ref{eqAdjEqsFirstCaseAdjDaeDiff})--(\ref{eqAdjEqsFirstCaseAdjDaeAlg}) is provided in the second part of the paper \cite{ps19b}.

\section{Numerical calculation of reduced gradients} \label{SecDiscrete}

In this chapter we present our approach to the numerical calculation of the reduced gradients of some endpoint state functional against the controls and discuss the convergence of the proposed numerical scheme. We focus on the numerical integration of continuous equations, as the numerical calculations at transition and end points reduces to solving the system of algebraic equations. The calculation of reduced gradients requires numerical integration of system equations, adjoint equations and reduced gradient formula. We discuss separately the ODEs case and DAEs case (which is presented in \cite{ps19b}). 

We start our considerations from a control system described by ODEs
\begin{equation}
x'(t) = f(x(t),u(t)),\ x(t_0) = x_0.
\label{discrEq:sysODE}
\end{equation}
The system equations (\ref{discrEq:sysODE}) are numerically integrated over the time interval $ [t_0, t_f] $ using the Runge-Kutta scheme (\cite{hnw93})
\begin{equation}
x_i(k+1) = x(k)+h(k)\sum_{j=1}^{s} a_{ij} f\left(x_j(k+1),\ u(t(k)+c_ih(k)) \right)
\label{discrEq:rkSchemeIntODE}
\end{equation}
for $ i = 1,\dots,s $ and
\begin{equation}
x(k+1) = x(k) + h(k)\sum_{i=1}^{s} b_i f\left(x_i(k+1),\ u(t(k)+c_ih(k)) \right).
\label{discrEq:rkSchemeStepODE}
\end{equation}
where $ t(k) $ is time at the actual step, $ h(k) $ is the actual step size and $ x(k)\simeq x(t(k)) $ is the numerical approximation of the state at $ t(k) $. The constant coefficients $ a_{ij} $, $ b_i $, $ c_i $ for $ i,j = 1,\dots,s $ define the Runge-Kutta scheme. At each step of the Runge-Kutta scheme the nonlinear system (\ref{discrEq:rkSchemeIntODE}) is first solved for variables $ x_i(k+1),\ i = 1,\dots,s $ and then (\ref{discrEq:rkSchemeStepODE}) is used to calculate $ x(k+1) $. The scheme is repeated for steps $ k = 0,\dots, K-1 $, where $ x(0) = x_0 $, $ t(0)=t_0 $ and $ t(K)=t_f $. In our code we assume the control variable is constant along the single integration step, so we have
\begin{equation}
u(t(k)+c_ih(k)) = u(k),\ i=1,\ldots,s.
\end{equation}

In the next step we must calculate numerically the adjoint equations %to (\ref{discrEq:sysODE}) are
\begin{equation}
\lambda'(t) = -f_x^T (x(t),u(t))\lambda(t),%\ \lambda(t_f) = \phi_x^T(x(t_f)).
\label{discrEqAdjODE}
\end{equation}
The adjoint equations are solved backwards in time with an endpoint condition given at time $ t_f $. The adjoint equations are ODEs, so the Runge-Kutta scheme can be used to numerically integrate them. The Runge-Kutta scheme for adjoint equations is
\begin{eqnarray}
\lambda_i(k) &=& \lambda(k+1) - \bar{h}(k+1) \sum_{j=1}^{\bar{s}} \bar{a}_{ij} \left(  -f_x^T 
\left( 
	x\left( \bar{t}(k+1)-\bar{c}_j\bar{h}(k+1) \right),
	\right . \right . \nonumber \\
& & \left . \left . \bar{u}(k+1) 
\right) \lambda_j(k) \right)
\nonumber\\
&=& \lambda(k+1) + \bar{h}(k+1) \sum_{j=1}^{\bar{s}} \bar{a}_{ij} f_x^T 
\left( 
x\left( \bar{t}(k+1)-\bar{c}_j\bar{h}(k+1) \right), \right .
\nonumber \\
& & \left .  \bar{u}(k+1) 
\right) \lambda_j(k) \nonumber \\
\label{discrEq:rkSchemeIntAdjODE} 
\end{eqnarray}
for $ i = 1,\dots,\bar{s} $ and
\begin{eqnarray}
\lambda(k) &=& \lambda(k+1) - \bar{h}(k+1) \sum_{i=1}^{\bar{s}} \bar{b}_i
\left( 
-f_x^T 
\left( 
	x\left( \bar{t}(k+1)-\bar{c}_i\bar{h}(k+1) \right), \right . \right .
	\nonumber \\
	& & \left . \left .  \bar{u}(k+1) 
\right) \lambda_i(k) \right) \nonumber \\
&=& \lambda(k+1) + \bar{h}(k+1) \sum_{i=1}^{\bar{s}} \bar{b}_i
f_x^T 
\left( 
	x\left( \bar{t}(k+1)-\bar{c}_i\bar{h}(k+1) \right),\right . \nonumber \\
& & \left .  \bar{u}(k+1) 
\right) \lambda_i(k), \nonumber \\
\label{discrEq:rkSchemeStepAdjODE} 
\end{eqnarray}
where $ \bar{t}(k) $ is time at the step $ k $, $ \bar{h}(k)>0 $ is the step size at step $ k $ and $ \lambda(k) \simeq \lambda(t(k)) $ is the approximation of the adjoint state at $ t(k) $. At each step of the Runge-Kutta scheme, the system of equations (\ref{discrEq:rkSchemeIntAdjODE}) is first solved for variables $ \lambda_i(k),\ i = 1,\dots,\bar{s} $ and then (\ref{discrEq:rkSchemeStepAdjODE}) is used to calculate $ \lambda(k) $. The scheme is repeated for steps $ k = \bar{K}-1,\dots, 0 $, where $ t(0)=t_0 $ and $ t(\bar{K})=t_f $. In general the number of steps, time moments at particular steps, step sizes and control values may be different for numerical integration of system and adjoint equations. It has been denoted by adding the bar symbol over the appropriate variables ($ \bar{K} $, $ \bar{t}(k) $, $ \bar{h}(k) $, $ \bar{u}(k) $) in (\ref{discrEq:rkSchemeIntAdjODE})-(\ref{discrEq:rkSchemeStepAdjODE}). What is more, the Runge-Kutta schemes used for numerical integration of system and adjoint equations can also be different, what has been emphasized by introducing the coefficients  $ \bar{a}_{ij} $, $ \bar{b}_i $, $ \bar{c}_i $ for $ i,j = 1,\dots,\bar{s} $ in the Runge-Kutta scheme for adjoint equations. The additional difficulty in implementation of the scheme (\ref{discrEq:rkSchemeIntAdjODE})-(\ref{discrEq:rkSchemeStepAdjODE}), is that a subprocedure for the reliable approximation of the state $ x(t) $ for $ t \in \left[ \bar{t}(k), \bar{t}(k+1) \right]  $ is required, such that terms $ x\left( \bar{t}(k+1)-\bar{c}_i\bar{h}(k+1) \right) $ are properly approximated. 

To derive the numerical scheme for calculation of reduced gradients, let us consider the differential of the endpoint state functional (corresponding to a perturbation $d$) of the form 
\begin{equation}
d \phi(x(t_f)) = \int_{t_0}^{t_f} d(t)^T f_u^T(x(t),u(t)) \lambda(t) dt.
\label{discrEq:funcDiffCont}
\end{equation}
The control signal is approximated as a piecewise constant function
\begin{equation}
u(t) = u_n\ for\ t\in\left( t_{n-1}, t_{n}\right],\ n = 1,\dots,N,
\end{equation}
where
\begin{equation}
t_{n} = t_0 + n\frac{t_f-t_0}{N}.
\end{equation}
Now (\ref{discrEq:funcDiffCont}) can be written as 
\begin{equation}
d \phi(x(t_f)) = \sum_{l = 1}^{N}  d^T_n \int_{t_{n-1}}^{t_{n}}  f_u^T(x(t),u_n) \lambda(t)  dt,
\label{discrEq:funcDiffContPWCControl}
\end{equation} 
where $ d_n $ is a differential of the control signal $ u_n $. The derivative of $ \phi(x(t_f)) $ with respect to $ u_n $ is therefore
\begin{equation}
\frac{d \phi(x(t_f))}{ du_n } = \int_{t_{n-1}}^{t_{n}} f_u^T(x(t),u_n) \lambda(t)  dt.
\label{discrEq:redGradCont}
\end{equation}
To approximate numerically the above integral, we choose times
\begin{equation}
\tilde{t}(0)<\dots<\tilde{t}(k)<\dots<\tilde{t}(\tilde{K}),
\end{equation}
such that $ \tilde{t}(0) = t_0 $, $ \tilde{t}(\tilde{K}) = t_f $ and for each $ n $ there is an index $ k_n $, such that 
\begin{equation}
\tilde{t}(k_n) = t_n.
\end{equation}  
The values $\tilde{t}(0),\dots,\tilde{t}(k) ,\dots,\tilde{t}(\tilde{K})  $ defines the division of the interval $ [t_0,t_f] $ into the subintervals. Now (\ref{discrEq:redGradCont}) can be written as
\begin{equation}
\frac{d \phi(x(t_f))}{ du_n } = \sum_{k=k_{n-1}}^{k_n-1} \int_{\tilde{t}(k)}^{\tilde{t}(k+1)} f_u^T(x(t),u_n) \lambda(t)  dt.
\label{discrEq:redGradContDiv}
\end{equation}
The integrals in (\ref{discrEq:redGradContDiv}) can be approximated by the quadrature formula
\begin{eqnarray}
&{\displaystyle \int_{\tilde{t}(k)}^{\tilde{t}(k+1)} f_u^T(x(t),u_n) \lambda(t)  dt \simeq 
\tilde{h}(k) \sum_{i=1}^{\tilde{s}} \tilde{b}_i f_u^T \left(x \left(\tilde{t}(k)+\tilde{c}_i \tilde{h}(k)\right),u_n \right)\times } \nonumber \\
&{\displaystyle  \lambda \left( \tilde{t}(k)+\tilde{c}_i \tilde{h}(k) \right),}
\label{discrEq:redGradQuadCont}
\end{eqnarray}
where $ \tilde{h}(k) = \tilde{t}(k+1)-\tilde{t}(k) $ is the step size at step $ k $ and $ \tilde{b}_i, \tilde{c}_i,\ i=1,\cdots,\tilde{s} $ are the quadrature coefficients. Again, the number of steps, time moments at particular steps and step sizes used for numerical integration of reduced gradients may be different from the respective values used for numerical integration of system or adjoint equations. It has been denoted by adding the tilde symbol over the appropriate variables ($ \tilde{K} $, $ \tilde{t}(k) $, $ \tilde{h}(k) $) in the above formulas. The additional difficulty in implementation of the quadrature formula (\ref{discrEq:redGradQuadCont}), is that a subprocedure for the reliable approximation of the state $ x(t) $ and adjoint state $ \lambda(t) $ for $ t \in \left[ \tilde{t}(k), \tilde{t}(k+1) \right]  $ is required, such that terms $ x\left( \tilde{t}(k)+\tilde{c}_i\tilde{h}(k) \right) $ and $ \lambda\left( \tilde{t}(k)+\tilde{c}_i\tilde{h}(k) \right) $ are properly approximated. 

In our code we utilize different approach to numerical calculation of adjoint equations and reduced gradients. The proposed approach is based on the concept of discrete state equations, discrete adjoint equations and discrete reduced gradients. Using that approach, the subprocedure for the approximation of the state $ x(t) $ and adjoint state $ \lambda(t) $ for $ t \in \left[ \tilde{t}(k), \tilde{t}(k+1) \right]  $ is not required. 
   
Let us start from formulation of the discrete equations. Let us rewrite the numerical integration scheme (\ref{discrEq:rkSchemeIntODE})-(\ref{discrEq:rkSchemeStepODE}) in a vector form
\begin{equation}
\left( 
\begin{array}{c}
	x_1 - x - h\sum_{j=1}^{s}a_{1j} f(x_j,u) \\
	\vdots \\
%	x_i - x - h\sum_{j=1}^{s}a_{ij} f(x_j,u)  \\
%	\vdots \\
	x_s - x - h\sum_{j=1}^{s}a_{sj} f(x_j,u) \\
	x^+ - x - h\sum_{i=1}^{s}b_i f(x_i,u) \\
\end{array}
\right) = 
\left( 
\begin{array}{c}
	0 \\ 
	\vdots \\ 
%	0 \\ 
%	\vdots \\ 
	0 \\ 
	0
\end{array}
\right).
\label{discrEq:rkSchemeVecODE}
\end{equation}
For the sake of a more compact notation we omitted the discrete step argument and introduced the symbols $ x_i = x_i(k+1) $, $ x = x(k) $, $ x^+ = x(k+1) $, $ h = h(k), u = u(k) $. If we now define the augmented state vector $ X(k) $ as
\begin{equation}
X(k) = \left( x_1(k)^T, \ldots ,x_s(k)^T, x(k)^T \right)^T ,\label{discrEq:VectorXODE}
\end{equation}
then (\ref{discrEq:rkSchemeVecODE}) can be presented in a form of the implicit discrete time state equation
\begin{equation}
F\left(X(k+1),X(k),u(k)\right) = 0 .
\label{discrEq:discrStateEq}
\end{equation}

Having the discrete state equation (\ref{discrEq:discrStateEq}) defined we may consider the discrete time adjoint equations (\cite{py11})
\begin{equation}
\Lambda(k) = -F_X^T(k)\left( F_{X^+}^{T} (k) \right) ^{-1}\Lambda(k+1) ,
\label{discrEq:discrAdjEq}
\end{equation}
where $ \Lambda(k) $ is the discrete adjoint variable at a discrete step $ k $ and
\begin{eqnarray}
F_{X^+}(k) &=& \frac{\partial F(X(k+1),X(k),u(k))}{\partial X(k+1)}, \\ 
F_{X}(k) &=& \frac{\partial F(X(k+1),X(k),u(k))}{\partial X(k)}  .
\end{eqnarray}
The equation (\ref{discrEq:discrAdjEq}) is equivalent to a system of equations
\begin{eqnarray}
F_{X^+}^{T} (k) R(k) &=& \Lambda(k+1) \label{discrEq:discrAdjEqRCalc}, \\
\Lambda(k) &=& -F_X^T(k) R(k) \label{discrEq:discrAdjEqLamCalc},
\end{eqnarray}  
where $ R(k+1) $ is an auxiliary variable introduced to avoid calculating the inverse of the matrix $ F_{X^+}^{T} (k) $.

Now we are going to rewrite the discrete adjoint equations in an element-wise manner. The partial derivatives matrices are 
\begin{equation}
F_{X^+}(k) = \left(
\begin{array}{cccc}
I - ha_{11}f_{x1} & \ldots & -ha_{1s}f_{xs} & 0  \\
\vdots &  & \vdots & \vdots  \\
-ha_{s1}f_{x1} & \ldots & I - ha_{ss}f_{xs} & 0  \\
-h b_1 f_{x1} & \ldots & -h b_s f_{xs}  & I 
\end{array}
\right ),\label{discrEq:FXPlusODE}
\end{equation}
\begin{equation}
F_{X}(k) = \left(
\begin{array}{cccc}
0 &  \ldots & 0 &  -I \\
\vdots & & \vdots & \vdots  \\
0 & \ldots & 0 & -I \\
0 & \ldots & 0 & -I  
\end{array}
\right ),\label{discrEq:FXODE}
\end{equation}
where for the sake of a shorter notation we denote 
\begin{eqnarray}
f_{xi} &=& f_x(x_i,u).
\end{eqnarray}
Let us denote 
\begin{eqnarray}
\Lambda(k) & = & \left( l_{1}(k)^T, \ldots, l_{s}(k)^T, \lambda(k)^T \right)^T,\nonumber \\
R(k) & = & \left( r_{1}(k)^T, \ldots, r_{s}(k)^T, r(k)^T \right)^T.\nonumber
\end{eqnarray}
Now (\ref{discrEq:discrAdjEqRCalc}) at a discrete time step $ k $  takes the form
\begin{equation}
\left(
\begin{array}{cccc}
	I - ha_{11}f^T_{x1} & \ldots & -ha_{s1}f^T_{x1} & -h b_1 f^T_{x1}   \\
	\vdots &  & \vdots & \vdots \\
	-ha_{1s}f^T_{xs} & \ldots & I - ha_{ss}f^T_{xs} & -h b_s f^T_{xs}   \\
	0 & \ldots & 0 & I  \\	
\end{array}
\right ) 
\left( 
\begin{array}{c}
	r_{1} \\ \vdots \\ 
	r_{s} \\  r 
\end{array}
\right) = 
\left(
\begin{array}{c}
	l^+_{1} \\  \vdots \\ 
	l^+_{s} \\  \lambda^+ 
\end{array}
\right) 
\label{discrEq:discrAdjEqRCalcVec}
\end{equation}
and (\ref{discrEq:discrAdjEqLamCalc}) takes the form
\begin{equation}
\left(
\begin{array}{c}
l_{1} \\ \vdots \\ 
l_{s} \\ \lambda  
\end{array}
\right) = -
\left(
\begin{array}{cccc}
0 & \ldots & 0 & 0  \\
\vdots & & \vdots & \vdots  \\
0 & \ldots & 0 & 0  \\
-I & \ldots & -I & -I  \\	
\end{array}
\right)
\left( 
\begin{array}{c}
r_{1} \\ \vdots \\ 
r_{s} \\ r
\end{array}
\right).
\label{discrEq:discrAdjEqLamCalcVec}
\end{equation}
For the sake of a more compact notation we omitted the discrete step argument and introduced the symbols $ r_{i} = r_{i}(k) $, $ r = r(k) $, $ l^+_{i} = l_{i}(k+1) $, $ \lambda^+ = \lambda(k+1) $, $ l_{i} = l_{i}(k) $, $ \lambda = \lambda(k) $.
If we rewrite (\ref{discrEq:discrAdjEqRCalcVec}) as a system of equations we obtain
\begin{equation}
r_{i} = \sum_{j=1}^{s} ha_{ji} f_{xi}^T r_{j} + hb_i f_{xi}^T r + l^+_{i} \label{discrEq:discrAdjEqRIntCalcSysEq001} 
\end{equation}
for $ i = 1,\ldots,s $ and
\begin{equation} 
r = \lambda^+. \label{discrEq:discrAdjEqREndCalcSysEq001} \\
\end{equation}
Using (\ref{discrEq:discrAdjEqREndCalcSysEq001}), (\ref{discrEq:discrAdjEqRIntCalcSysEq001}) can be written as
\begin{equation}
r_{i} = \sum_{j=1}^{s} ha_{ji} f_{xi}^T r_{j} + hb_i f_{xi}^T \lambda^+ + l^+_{i} . \label{discrEq:discrAdjEqRIntCalcSysEq002}
\end{equation}
If we rewrite (\ref{discrEq:discrAdjEqLamCalcVec}) as a system of equations we get
\begin{equation}
l_i = 0 \label{discrEq:discrAdjEqLamIntCalcSysEq001} 
\end{equation}
for $ i = 1,\ldots,s $ and
\begin{equation}
\lambda = \sum_{i=1}^{s} r_{i} + r. \label{discrEq:discrAdjEqLamEndCalcSysEq001} 
\end{equation}
Using (\ref{discrEq:discrAdjEqREndCalcSysEq001}), (\ref{discrEq:discrAdjEqLamEndCalcSysEq001}) becomes
\begin{equation}
\lambda = \sum_{i=1}^{s} r_{i} + \lambda^+. \label{discrEq:discrAdjEqLamEndCalcSysEq002} 
\end{equation}
From (\ref{discrEq:discrAdjEqLamIntCalcSysEq001}) we have that  $ l_{1}(k) = 0, \ldots,l_{s}(k) = 0$ for steps $ k = 0,...,K-1 $. During the analysis we also assume that $ l_{1}(K) = 0, \ldots,l_{s}(K) = 0$ . 
Under that assumption (\ref{discrEq:discrAdjEqRIntCalcSysEq002}) is equivalent to
\begin{equation}
r_{i} = \sum_{j=1}^{s} ha_{ji} f_{xi}^T r_{j} + hb_i f_{xi}^T \lambda^+  \label{discrEq:discrAdjEqRIntCalcSysEq003} 
\end{equation}
for $ i = 1,\dots,s $.

To show the resemblance of discrete adjoint equations and the Runge-Kutta scheme for continuous adjoint equations, we introduce the auxiliary variables in a manner proposed in (\cite{h00})
\begin{equation}
\lambda_{i} = \sum_{j=1}^{s} \frac{a_{ji}}{b_i}  r_{j} + \lambda^+ \label{discrEq:discrAdjEqRIntCalcSysEq004}   
\end{equation}
for $ i = 1,\ldots,s $. We have then 
\begin{equation}
r_{i} = hb_i f_{xi}^T\lambda_{i} \label{discrEq:discrAdjEqRIntCalcSysEq041}
\end{equation}
and (\ref{discrEq:discrAdjEqRIntCalcSysEq004}) can be transformed to the following form
\begin{equation}
\lambda_{i} =  h\sum_{j=1}^{s} \frac{a_{ji}b_j}{b_i} f_{xj}^T \lambda_{j} + \lambda^+ \label{discrEq:discrAdjEqRIntCalcSysEq005} 
\end{equation} 
and (\ref{discrEq:discrAdjEqLamEndCalcSysEq002}) can be transformed to
\begin{equation}
\lambda = h\sum_{i=1}^{s} b_i f_{xi}^T \lambda_i + \lambda^+ . \label{discrEq:discrAdjEqLamEndCalcSysEq003}
\end{equation} 
Let us rewrite (\ref{discrEq:discrAdjEqRIntCalcSysEq005}) and (\ref{discrEq:discrAdjEqLamEndCalcSysEq003}) by inserting the discrete step argument
\begin{equation}
\lambda_{i}(k) = \lambda(k+1) + h(k)\sum_{j=1}^{s} \frac{a_{ji}b_j}{b_i} f_{x}^T(x_j(k+1),u(k)) \lambda_{j}(k)
\label{discrEq:discrAdjEqRIntCalcSysEq006}  
\end{equation}
for $ i = 1,\dots,s $ and
\begin{equation}
\lambda(k) = \lambda(k+1) + h(k)\sum_{i=1}^{s} b_i f_{x}^T(x_i(k+1),u(k)) \lambda_i(k)  . \label{discrEq:discrAdjEqLamEndCalcSysEq004}
\end{equation}

For the convenience of reading we recall here the Runge-Kutta scheme for adjoint equations (\ref{discrEq:rkSchemeIntAdjODE})-(\ref{discrEq:rkSchemeStepAdjODE})
\begin{equation}
\lambda_i(k) = \lambda(k+1) + \bar{h}(k+1) \sum_{j=1}^{\bar{s}} \bar{a}_{ij} f_x^T 
\left( 
x\left( \bar{t}(k+1)-\bar{c}_j\bar{h}(k+1) \right), \bar{u}(k+1) 
\right) \lambda_j(k) \nonumber
\end{equation}
for $ i = 1,\dots,s $ and
\begin{equation}
\lambda(k) = \lambda(k+1) + \bar{h}(k+1) \sum_{i=1}^{\bar{s}} \bar{b}_i
f_x^T 
\left( 
x\left( \bar{t}(k+1)-\bar{c}_i\bar{h}(k+1) \right), \bar{u}(k+1) 
\right) \lambda_i(k). \nonumber
\end{equation}
Let us suppose, that the adjoint equations are integrated such that discrete times are the same as during the numerical integration of system equations so
\begin{equation}
\bar{t}(k) = t(k).
\end{equation} 
In such case we obtain
\begin{eqnarray}
\bar{h}(k+1) &=& h(k), \\
\bar{u}(k+1) &=& u(k)
\end{eqnarray}
and
\begin{equation}
\bar{t}(k+1)-\bar{c}_i\bar{h}(k+1) = t(k+1) - \bar{c}_ih(k) = t(k)+(1-\bar{c}_i)h(k).
\end{equation}
By introducing the above relations into (\ref{discrEq:rkSchemeIntAdjODE})-(\ref{discrEq:rkSchemeStepAdjODE}) we obtain
\begin{equation}
\lambda_i(k) = \lambda(k+1) + h(k) \sum_{j=1}^{\bar{s}} \bar{a}_{ij} f_x^T 
\left( 
	x\left( t(k)+(1-\bar{c}_j)h(k) \right), u(k) 
\right) \lambda_j(k) \label{discrEq:rkSchemeIntAdjODE001}
\end{equation}
for $ i = 1,\dots,s $ and
\begin{equation}
\lambda(k) = \lambda(k+1) + h(k) \sum_{i=1}^{\bar{s}} \bar{b}_i
f_x^T 
\left( 
	x\left( t(k)+(1-\bar{c}_i)h(k) \right), u(k) 
\right) \lambda_i(k). \label{discrEq:rkSchemeStepAdjODE001}
\end{equation}
If we now assume that
\begin{equation}
\bar{s} = s,\ \bar{a}_{ij} = \frac{a_{ji}b_j}{b_i},\ \bar{b}_i = b_i,\ \bar{c}_i = 1-c_i,\ i,j = 1,\dots,s \label{discrEq:rkSchemeCoeffAdj}
\end{equation}
(\ref{discrEq:rkSchemeIntAdjODE001})-(\ref{discrEq:rkSchemeStepAdjODE001}) can be rewritten as
\begin{equation}
\lambda_i(k) = \lambda(k+1) + h(k) \sum_{j=1}^{s} \frac{a_{ji}b_j}{b_i} f_x^T 
\left( 
	x\left( t(k)+c_jh(k) \right), u(k) 
\right) \lambda_j(k) \label{discrEq:rkSchemeIntAdjODE002}
\end{equation}
for $ i = 1,\dots,s $ and
\begin{equation}
\lambda(k) = \lambda(k+1) + h(k) \sum_{i=1}^{s} b_i
f_x^T \left( 
	x\left( t(k)+c_ih(k) \right), u(k) 
\right) \lambda_i(k). \label{discrEq:rkSchemeStepAdjODE002}
\end{equation}
Now the Runge-Kutta scheme (\ref{discrEq:rkSchemeIntAdjODE002})-(\ref{discrEq:rkSchemeStepAdjODE002}) for adjoint equations is almost identical to discrete adjoint equations (\ref{discrEq:discrAdjEqRIntCalcSysEq006})-(\ref{discrEq:discrAdjEqLamEndCalcSysEq004}). The only difference is that in discrete equations $ x_i(k+1) $ is used instead of $ x(t(k)+c_ih(k)) $. The convergence of (\ref{discrEq:rkSchemeIntAdjODE002})-(\ref{discrEq:rkSchemeStepAdjODE002}) is guaranteed if only the Runge-Kutta scheme $ \bar{a}_{ij} ,\ \bar{b}_i,\ \bar{c}_i,\ i,j = 1,\dots,s $ satisfies appropriate conditions. We will also justify that the usage of $ x_i(k+1) $ instead of $ x(t(k)+c_ih(k)) $ does not break the convergence of (\ref{discrEq:discrAdjEqRIntCalcSysEq006})-(\ref{discrEq:discrAdjEqLamEndCalcSysEq004}) to continuous adjoint trajectory. 

The discrete reduced gradients are defined as (\cite{py11})
\begin{equation}
\frac{d \phi(X(K))}{du_n} =  \sum_{k=k_{n-1}}^{k_n-1} -F_u^T(k)\left( F_{X^+}^{T} (k) \right) ^{-1}\Lambda(k+1) = \sum_{k=k_{n-1}}^{k_n-1} -F_u^T(k)R(k+1),
\label{discrEq:discrRedGrad}
\end{equation}
where
\begin{equation}
F_{u}(k) = \frac{\partial F(X(k+1),X(k),u(k))}{\partial u(k)}.
\end{equation}
The above partial derivative rewritten in the element-wise manner is
\begin{equation}
F_{u}(k) = \left(
\begin{array}{c}
- h\sum_{j=1}^{s}a_{1j} f_{uj}  \\
\vdots  \\
- h\sum_{j=1}^{s}a_{sj} f_{uj} \\
- h\sum_{i=1}^{s}b_{i} f_{ui}
\end{array}
\right ),
\end{equation}
where for the sake of the shorter notation we omitted the discrete step argument and introduced
\begin{equation}
f_{ui} = f_u(x_i,u).
\end{equation}
Let us now inspect $ -F_u^T(k)R(k+1) $  %Now (\ref{discrEq:discrRedGrad}) is
\begin{eqnarray}
-F_u^T(k)R(k+1) &=&  
- \left(
\begin{array}{c}
- h\sum_{j=1}^{s}a_{1j} f_{uj}  \\
\vdots  \\
- h\sum_{j=1}^{s}a_{sj} f_{uj} \\
- h\sum_{i=1}^{s}b_{i} f_{ui}
\end{array}
\right )^T
\left( 
\begin{array}{c}
r_{1} \\ \vdots \\ 
r_{s} \\ r
\end{array}
\right) \\
&=&  \sum_{i = 1}^{s}  h \left(\sum_{j=1}^{s}a_{ij}  f^T_{uj} \right) r_i + h\sum_{i=1}^{s}b_{i} f^T_{ui} r. \nonumber
\end{eqnarray}
From (\ref{discrEq:discrAdjEqREndCalcSysEq001}) and (\ref{discrEq:discrAdjEqRIntCalcSysEq041}) we get
\begin{eqnarray}
-F_u^T(k)R(k+1) &=&  \sum_{i = 1}^{s}  h \left(\sum_{j=1}^{s}a_{ij}  f^T_{uj} \right) hb_i f_{xi}^T\lambda_{i} + h\sum_{i=1}^{s}b_{i} f^T_{ui} \lambda^+ \\
&=& h \sum_{j = 1}^{s} f^T_{uj} \left( h \sum_{i=1}^{s}a_{ij} b_i f_{xi}^T\lambda_{i} \right)  + h\sum_{i=1}^{s}b_{i} f^T_{ui} \lambda^+ \nonumber \\
&=& h \sum_{j = 1}^{s} f^T_{uj}b_j \left( h \sum_{i=1}^{s} \frac{a_{ij} b_i}{b_j} f_{xi}^T\lambda_{i} \right)  + h\sum_{i=1}^{s}b_{i} f^T_{ui} \lambda^+ . \nonumber
\end{eqnarray}
From (\ref{discrEq:discrAdjEqRIntCalcSysEq005}) we get that
\begin{eqnarray}
-F_u^T(k)R(k+1) &=&  h \sum_{j = 1}^{s} f^T_{uj}b_j \left( \lambda_j - \lambda^+ \right)  + h\sum_{i=1}^{s}b_{i} f^T_{ui} \lambda^+ \\
&=&  h \sum_{i = 1}^{s} f^T_{ui}b_i \left(  \lambda_i - \lambda^+ \right)  + h\sum_{i=1}^{s}b_{i} f^T_{ui} \lambda^+ \nonumber\\
&=&  h \sum_{i = 1}^{s} b_i f^T_{ui}  \lambda_i . \nonumber
\end{eqnarray} 
The discrete reduced gradient now takes the form
\begin{equation}
\frac{d \phi(X(K))}{du_n} =  \sum_{k=k_{n-1}}^{k_n-1} h(k) \sum_{i = 1}^{s} b_i f_u^T(x_i(k+1),u_n)  \lambda_i(k).
\label{discrEq:discrRedGrad002}
\end{equation}
The quadrature formula for the continuous reduced gradient calculation is
\begin{eqnarray}
\left( \frac{d \phi(x(t_f))}{ du_n } \right)_Q & = & \sum_{k=k_{n-1}}^{k_n-1}
\tilde{h}(k) \sum_{i=1}^{\tilde{s}} \tilde{b}_i f_u^T \left(x \left(\tilde{t}(k)+\tilde{c}_i \tilde{h}(k)\right),u_n \right)\times \nonumber \\ 
& & \lambda \left( \tilde{t}(k)+\tilde{c}_i \tilde{h}(k) \right).
\label{discrEq:redGradQuadCont002}
\end{eqnarray}
Let us now suppose, that the quadrature is calculated such that the discrete times are the same as during the numerical integration of system equations so
\begin{equation}
\tilde{t}(k) = t(k).
\end{equation} 
In such case we obtain $ \tilde{h}(k) = h(k) $ and
\begin{equation}
\left( \frac{d \phi(x(t_f))}{ du_n } \right)_Q = \sum_{k=k_{n-1}}^{k_n-1}
h(k) \sum_{i=1}^{\tilde{s}} \tilde{b}_i f_u^T \left(x \left(t(k)+\tilde{c}_i h(k)\right),u_n \right) \lambda \left( t(k)+\tilde{c}_i h(k) \right).
\label{discrEq:redGradQuadCont003}
\end{equation}
Let us also assume that
\begin{equation}
\tilde{s} = s,\ \tilde{b}_i = b_i,\ \tilde{c}_i = c_i,\ i = 1,\dots,s. \label{discrEq:quadCoeffDiscrRedGrad}
\end{equation}
In such case (\ref{discrEq:redGradQuadCont003}) becomes
\begin{equation}
\left( \frac{d \phi(x(t_f))}{ du_n } \right)_Q = \sum_{k=k_{n-1}}^{k_n-1}
h(k) \sum_{i=1}^{s} b_i f_u^T \left(x \left(t(k)+c_i h(k)\right),u_n \right) \lambda \left( t(k)+c_i h(k) \right).
\label{discrEq:redGradQuadCont004}
\end{equation}
Now the quadrature scheme (\ref{discrEq:redGradQuadCont004}) for continuous reduced gradients is almost identical to discrete reduced gradients  (\ref{discrEq:discrRedGrad002}). The only difference is that in (\ref{discrEq:discrRedGrad002}) $ x_i(k+1) $ and $ \lambda_i(k) $ are used instead of $ x(t(k)+c_ih(k)) $ and $ \lambda(t(k)+c_ih(k)) $. The convergence of (\ref{discrEq:redGradQuadCont004}) is guaranteed if only the quadrature with coefficients $ b_i,\ c_i,\ i = 1,\dots,s $ satisfies appropriate conditions. We will also justify that the usage of  $ x_i(k+1) $ and $ \lambda_i(k) $ instead of $ x(t(k)+c_ih(k)) $ and $ \lambda(t(k)+c_ih(k)) $ does not break the convergence of (\ref{discrEq:discrRedGrad002}) to continuous reduced gradients. 

After showing the correspondence between discrete equations and their continuous counterparts, we want to justify the convergence of the discrete equations. To obtain the right error order estimated we always assume that the system function $ f(x,u) $ and its partial derivatives are Lipschitz continuous functions. The convergence analysis for any Runge-Kutta scheme is a challenging task, so we limit ourselves to the analysis of the RADAU IIA integration scheme for $s=3$. This scheme has been chosen for the numerical integration of the system equations because of its advantageous properties, especially for the integration of DAEs (\cite{hlr89}).
The coefficients of the RADAU IIA scheme for $ s=3 $ collected in matrices are (\cite{hw96}, p. 74)
\begin{equation}
A = (a_{ij}) = 
\left(
\begin{array}{ccc}
	\frac{11}{45}-\frac{7\sqrt{6}}{360} & \frac{37}{225}-\frac{169\sqrt{6}}{1800} & -\frac{2}{225}+\frac{\sqrt{6}}{75} \\
	 & & \\
	\frac{37}{225}+\frac{169\sqrt{6}}{1800} & \frac{11}{45}+\frac{7\sqrt{6}}{360} & -\frac{2}{225}-\frac{\sqrt{6}}{75}\\
	 & & \\
	\frac{4}{9}-\frac{\sqrt{6}}{36} & \frac{4}{9}+\frac{\sqrt{6}}{36} & \frac{1}{9}
\end{array}
\right),
\end{equation}
\begin{equation*}
b = (b_i) = \left(
\begin{array}{ccc}
	\frac{4}{9}-\frac{\sqrt{6}}{36}, & \frac{4}{9}+\frac{\sqrt{6}}{36}, & \frac{1}{9}
\end{array} \right),\
c = (c_i)= \left(
\begin{array}{ccc}
	\frac{2}{5}-\frac{\sqrt{6}}{10}, & \frac{2}{5}+\frac{\sqrt{6}}{10}, & 1
\end{array} \right).
\end{equation*}
To calculate the convergence order $ p $ of the Runge-Kutta scheme the following conditions imposed on the Runge-Kutta coefficients are checked (\cite{hnw93}, p. 208) 
\begin{eqnarray}
B(p)&:&\ \forall_{l = 1,\dots,p}\ \sum_{i=1}^{s} b_ic_i^{l-1} = \frac{1}{l}, \\
C(q)&:&\ \forall_{l = 1,\dots,q}\ \forall_{i = 1,\dots,s}\ \sum_{j=1}^{s} a_{ij} c_j^{l-1} = \frac{c_i^l}{l}, \\
D(r)&:&\ \forall_{l = 1,\dots,r}\ \forall_{j = 1,\dots,s}\ \sum_{i=1}^{s} b_i c_i^{l-1} a_{ij}= \frac{b_j}{l} \left(1 - c_j^l\right),  
\end{eqnarray}
If the conditions $ B(p), C(q), D(r) $ are satisfied for $ p\leq 2q+2 $ and $ p\leq q+r+1 $ then the order of convergence of the Runge-Kutta scheme is $ p $ (\cite{hnw93}, Theorem 7.4 p. 208), i.e. the difference between approximated state $ x(k) $ and exact state at time $ t(k) $ satisfies
\begin{equation}
x(k) - x(t(k)) = O(h_m^p), \label{discrEq:errEstStateODEStep} 
\end{equation}
where 
\begin{equation}
h_m = max \{ h(0),\dots, h(K-1) \}.
\end{equation}
The RADAU IIA scheme satisfies conditions $ B(p=5), C(q=3), D(r=2) $ (\cite{hw96}, p. 77) and so the order of the method is $ p=5 $.

To justify the convergence of the discrete equations we consider the unperturbed (nominal) Runge-Kutta scheme 
\begin{eqnarray}
x_i^{n}(k+1) &=& x^n(k) + h(k)\sum_{j=1}^{s} a_{ij} f\left(x_j^{n}(k+1), u(k) \right),\ i = 1,\dots,s \label{discrEq:rkSchemeIntODENom} \\
x^{n}(k+1) &=& x^n(k) + h(k)\sum_{i=1}^{s} b_i f\left(x_i^{n}(k), u(k) \right) \label{discrEq:rkSchemeStepODENom}
\end{eqnarray}
and the perturbed Runge-Kutta scheme
\begin{eqnarray}
x_i^{p}(k+1) &=& x^n(k) + h(k)\sum_{j=1}^{s} a_{ij} f\left(x_j^{p}(k+1), u(k) \right) + \delta_i,\ i = 1,\dots,s \label{discrEq:rkSchemeIntODEPert} \\
x^{p}(k+1) &=& x^n(k) + h(k)\sum_{i=1}^{s} b_i f\left(x_i^{p}(k), u(k) \right) + \delta_{s+1} . \label{discrEq:rkSchemeStepODEPert}
\end{eqnarray}
where $ \delta_i,\ i=1,\dots,s+1 $ are perturbations. The estimates of the errors in the perturbed system are (\cite{hw96}, Theorem 14.3 p. 219)
\begin{eqnarray}
 x^p_i(k+1) - x^n_i(k+1) &=& O(\delta) ,\ i=1,\dots,s \label{discrEq:errEstStateODEIntPert} \\
 x^{p}(k+1) - x^{n}(k+1) &=& O(\delta) , \label{discrEq:errEstStateODEStepPert}  
\end{eqnarray}
where $ \delta = max\{ \delta_1,\dots, \delta_{s+1}\} $. It should be noted that the error estimates has been formulated for stiff differential equations (\cite{hw96}, Theorem 14.3 p. 219). For that reason some additional conditions has been imposed on the Runge-Kutta scheme. Those conditions concern the invertibility and coercivity of the matrix $ A=(a_{ij}) $. It is possible that those conditions could be relaxed for non-stiff systems, but as the considered schemes (RADAU IA and RADAU IIA) satisfy them (\cite{hw96}, Theorems 14.3 and 14.5) we do not discuss this issue further. 
The error estimates (\ref{discrEq:errEstStateODEIntPert})-(\ref{discrEq:errEstStateODEStepPert}) play a crucial role for the proposed convergence analysis. To use these estimates nominal and perturbed schemes must be formulated in an appropriate way in each case. 

Let us now derive the estimate of the global error $  x_i(k+1) - x(t(k) + c_ih(k))  $. The nominal Runge-Kutta scheme is formulated assuming that the exact state value is known at $ t(k) $ 
\begin{eqnarray}
x^n_i(k+1) &=& x(t(k)) + h(k)\sum_{j=1}^{s} a_{ij} f\left(x^n_j(k+1),\ u(k) \right),\ i = 1,\dots,s.\nonumber \\
& &  \label{discrEq:rkSchemeIntODENom001}  
\end{eqnarray}
If $ C(q) $ holds, the following local error estimate is valid (\cite{hnw93}, Lemma 7.5 p. 210) 
\begin{equation}
 x^n_i(k+1) - x(t(k) + c_ih(k))  = O(h^{q+1}(k)),\ i = 1,\dots,s. \label{discrEq:errEstStateODEIntLoc}
\end{equation}
In this case the perturbed Runge-Kutta scheme is actually the regular Runge-Kutta scheme, with the approximation of the state $ x(k) $ at time $ t(k) $ used
\begin{eqnarray}
x_i(k+1) &=& x(k) + h(k)\sum_{j=1}^{s} a_{ij} f\left(x_j(k+1),\ u(k) \right)  \label{discrEq:rkSchemeIntODEPert001} \\
&=& x(t(k)) + h(k)\sum_{j=1}^{s} a_{ij} f\left(x_j(k+1),\ u(k) \right) + (x(k)- x(t(k))) \nonumber 
\end{eqnarray}
for $ i = 1,\dots,s $, so the perturbations are
\begin{equation}
\delta_i = \delta = x(k)- x(t(k)),\ i = 1,\dots,s.
\end{equation}
From (\ref{discrEq:errEstStateODEStep}) we have $  \delta  =  x(k)- x(t(k))  =  O(h_m^p) $. From (\ref{discrEq:errEstStateODEIntPert}) we obtain
\begin{equation}
 x_i(k+1) - x^n_i(k+1) = O(h_m^p). \label{discrEq:errEstStateODEIntPert001}
\end{equation} 
%for some positive constant $ C_2 $.
By combining (\ref{discrEq:errEstStateODEIntLoc}) and (\ref{discrEq:errEstStateODEIntPert001})  we obtain the required global error estimate
\begin{eqnarray}
x_i(k+1) - x(t(k) + c_ih(k)) &=& O(h^{q+1}(k))+ O(h_m^p) = O(h^{q+1}_m)+ O(h_m^p) \nonumber \\
&=& O(h_m^{min\{p,q+1\}}).\label{discrEq:errEstXIntGlob}
\end{eqnarray}
For the RADAU IIA scheme with $ s=3 $ we get
\begin{equation}
x_i(k+1) - x(t(k) + c_ih(k)) = O(h_m^{min\{5,3+1\}}) = O(h_m^4).
\label{discrEq:errEstStateODEIntGlob}
\end{equation}

Let us now consider the Runge--Kutta scheme for adjoint equations (\ref{discrEq:rkSchemeIntAdjODE002})-(\ref{discrEq:rkSchemeStepAdjODE002}). The coefficients $ \bar{a}_{ij} ,\ \bar{b}_i,\ \bar{c}_i,\ i,j = 1,\dots,s $ of the Runge-Kutta scheme collected in matrices are 
\begin{equation}
\bar{A} = (\bar{a}_{ij}) = \left(
\begin{array}{ccc}
\frac{11}{45}-\frac{7\sqrt{6}}{360} & \frac{11}{45}+\frac{43\sqrt{6}}{360} & \frac{1}{9} \\
 & & \\
\frac{11}{45}-\frac{43\sqrt{6}}{360} & \frac{11}{45}-\frac{7\sqrt{6}}{360} & \frac{1}{9} \\
 & & \\
-\frac{1}{18}+\frac{\sqrt{6}}{18} & -\frac{1}{18}-\frac{\sqrt{6}}{18} & \frac{1}{9}
\end{array} \right),
\end{equation}
\begin{equation*}
\bar{b} = (\bar{b}_i) = \left(
\begin{array}{ccc}
\frac{4}{9}-\frac{\sqrt{6}}{36}, & \frac{4}{9}+\frac{\sqrt{6}}{36}, & \frac{1}{9}
\end{array} \right),\
\bar{c} = (\bar{c}_i) = \left(
\begin{array}{ccc}
\frac{3}{5}+\frac{\sqrt{6}}{10}, & \frac{3}{5}-\frac{\sqrt{6}}{10}, & 0
\end{array} \right). 
\end{equation*}
It turns out that the above Runge-Kutta scheme is RADAU IA scheme  (\cite{hw96}, p. 73). The RADAU IA scheme satisfies conditions $ B(\bar{p}=5), C(\bar{q}=2), D(\bar{r}=3) $ (\cite{hw96}, p. 77) and so the order of the method is $ \bar{p}=5 $. 

Let us now define the unperturbed Runge-Kutta scheme as the Runge-Kutta scheme for continuous adjoint equations with the exact adjoint state known at time $ t(k+1) $
\begin{equation}
\lambda_i^n(k) = \lambda(t(k+1)) + h(k) \sum_{j=1}^{s} \bar{a}_{ij}  f_x^T 
\left( 
x\left( t(k)+c_jh(k) \right), u(k) 
\right) \lambda_j^n(k) \label{discrEq:rkSchemeIntAdjODELoc}
\end{equation}
for $ i = 1,\dots,s $ and
\begin{equation}
\lambda^n(k) = \lambda(t(k+1)) + h(k) \sum_{i=1}^{s} \bar{b}_i
f_x^T \left( 
x\left( t(k)+c_ih(k) \right), u(k) 
\right) \lambda_i^n(k). \label{discrEq:rkSchemeStepAdjODELoc}
\end{equation}
As $ \bar{q}=2 $ and $ \bar{p}=5 $ the following local error estimates are valid
\begin{eqnarray}
\lambda^n_i(k) - \lambda(t(k)+c_ih(k)) &=& O(h^{\bar{q}+1}(k)),\ i=1,\dots,s\ ,
\label{discrEq:errEstAdjIntLoc} \\
\lambda^n(k) - \lambda(t(k)) &=& O(h^{\bar{p}+1}(k)).
\label{discrEq:errEstAdjEndLoc}
\end{eqnarray}
The perturbed Runge-Kutta scheme is defined as the discrete adjoint equations with the exact adjoint state known at time $ t(k+1) $ 
\begin{equation}
\lambda_{i}^p(k) = \lambda(t(k+1)) + h(k)\sum_{j=1}^{s} \bar{a}_{ij} f_{x}^T(x_j(k+1),u(k)) \lambda_{j}^p(k)
\label{discrEq:discrAdjEqRIntLoc}  
\end{equation}
for $ i = 1,\dots,s $ and
\begin{equation}
\lambda^p(k) = \lambda(t(k+1)) + h(k)\sum_{i=1}^{s} \bar{b}_{i} f_{x}^T(x_i(k+1),u(k)) \lambda_i^p(k)  . \label{discrEq:discrAdjEqEndLoc}
\end{equation}
Using (\ref{discrEq:errEstXIntGlob}), (\ref{discrEq:discrAdjEqRIntLoc})-(\ref{discrEq:discrAdjEqEndLoc}) can be rewritten as
\begin{equation}
\lambda_i^p(k) = \lambda(t(k+1)) + h(k) \sum_{j=1}^{s} \bar{a}_{ij} f_x^T 
\left( 
x\left( t(k)+c_jh(k) \right), u(k) 
\right) \lambda_j^p(k) + \delta_i \label{discrEq:discrAdjEqRIntLoc001}
\end{equation}
for $ i = 1,\dots,s $ and
\begin{equation}
\lambda^p(k) = \lambda(t(k+1)) + h(k) \sum_{i=1}^{s} \bar{b}_{i}
f_x^T \left( 
x\left( t(k)+c_ih(k) \right), u(k) 
\right) \lambda_i^p(k) + \delta_{s+1}, \label{discrEq:discrAdjEqEndLoc001}
\end{equation}
where the perturbations result from differences $ x_i(k+1) - x\left( t(k)+c_ih(k) \right) $ and so the perturbations satisfy
\begin{equation}
\delta_i = O(h(k)h^{q+1}_m),\ i=1,\dots,s+1.
\end{equation}
From (\ref{discrEq:errEstStateODEStepPert}) we obtain 
\begin{equation}
\lambda^p(k) - \lambda^n(k) = O(h(k)h^{q+1}_m).
\label{discrEq:errEstDiscrAdjEndPert}
\end{equation}
By combining (\ref{discrEq:errEstAdjEndLoc}) and (\ref{discrEq:errEstDiscrAdjEndPert}) we obtain the following local error estimate for discrete adjoint equations
\begin{equation}
\lambda^p(k) - \lambda(t(k)) = O(h^{\bar{p}+1}(k)) + O(h(k)h^{q+1}_m) = O(h(k)h^{min\{ \bar{p} ,q+1\}}_m).
\label{discrEq:errEstDiscrAdjEndLoc}
\end{equation} 
If at each step the local error of the numerical integration scheme is $ O(h(k)h^{p}_m) $, then the global error of the scheme is $ O(h^p_m) $ (\cite{hnw93}, Theorem 3.6 p. 162).   
Having the local error estimate (\ref{discrEq:errEstDiscrAdjEndLoc}) the global error estimate for the discrete adjoint equations (\ref{discrEq:discrAdjEqRIntCalcSysEq006})-(\ref{discrEq:discrAdjEqLamEndCalcSysEq004}) is therefore
\begin{equation}
\lambda(k) - \lambda(t(k)) = O(h^{min\{ \bar{p} ,q+1\}}_m) =  O(h^{\bar{p}_d}_m),
\label{discrEq:errEstDiscrAdjEndGlob}
\end{equation} 
where $ \bar{p}_d = min\{ \bar{p}, q+1 \} $. For RADAU IIA scheme we obtain $ \bar{p}_d = min\{ 5, 3+1 \} = 4$.
Let us now derive the global error estimate of $ \lambda_i(k) - \lambda(t(k)+c_ih(k)) $. The unperturbed system is  (\ref{discrEq:rkSchemeIntAdjODELoc}) and the perturbed system is (\ref{discrEq:discrAdjEqRIntCalcSysEq006}). The perturbed system can be rewritten as
\begin{equation}
\lambda_i(k) = \lambda(k) + h(k) \sum_{j=1}^{s} \bar{a}_{ij}  f_x^T 
\left( 
x_i(k+1), u(k) 
\right) \lambda_j(k) + \delta_i \label{discrEq:discrAdjEqIntPert}
\end{equation}
for $ i = 1,\dots,s $. The perturbations result from differences $ \lambda(k) - \lambda(t(k+1)) $  and $ x_i(k+1) - x\left( t(k)+c_ih(k) \right) $ and so the perturbations satisfy
\begin{eqnarray}
\delta_i = O(h^{\bar{p}_d}_m) + O(h(k)h^{q+1}_m) = O(h^{\bar{p}_d}_m) + O(h^{q+2}_m) = O(h_m^{min\{ \bar{p}_d, q+2 \}}).
\end{eqnarray}
From (\ref{discrEq:errEstStateODEIntPert}) we obtain
\begin{equation}
\lambda_i(k) - \lambda_i^n(k) = O(h_m^{min\{ \bar{p}_d, q+2 \}}).
\label{discrEq:errEstDiscrAdjIntPert}
\end{equation}
By combining (\ref{discrEq:errEstAdjIntLoc}) and (\ref{discrEq:errEstDiscrAdjIntPert}) we obtain
\begin{eqnarray}
\lambda_i(k) - \lambda(t(k)+c_ih(k)) &=& O(h^{\bar{q}+1}(k)) + O(h_m^{min\{ \bar{p}_d, q+2 \}}) \label{discrEq:errEstAdjIntGlob}\\
 &=& O(h_m^{min\{\bar{q}+1, \bar{p}_d, q+2 \}}) = O(h_m^{\bar{q}_d+1}), \nonumber
\end{eqnarray}
where $ \bar{q}_d = min\{\bar{q}+1, \bar{p}_d, q+2 \} -1 $. For RADAU IIA scheme we obtain $ \bar{q}_d = min\{2+1, 4, 3+2 \} -1 = 2 $.

Let us finally derive the order of the discrete reduced gradients (\ref{discrEq:discrRedGrad002}). If the quadrature coefficients $ \tilde{b}_i, \tilde{c}_i,\ i=1,\dots,s $ satisfy the condition $ B(\tilde{p}) $, then the quadrature integrates exactly the polynomials of a degree up to $ \tilde{p}-1 $ on an interval $ [0,1] $ (\cite{hlr89}, p.15). We consider the quadrature with coefficients $ \tilde{b}_i = b_i, \tilde{c}_i = c_i,\ i=1,\dots,s $, which satisfy the condition $ B(\tilde{p} = 5) $. The local error resulting from the usage of the quadrature (\ref{discrEq:redGradQuadCont004}) over single integration step can therefore be estimated as  
\begin{eqnarray}
&\int_{t(k)}^{t(k+1)}& f_u^T(x(t),u_n) \lambda(t)  dt \\
&=& \int_{0}^{1} f_u^T(x(t(k)+\tau h(k),u_n) \lambda(t(k)+\tau h(k)) h(k) d\tau \nonumber \\
&=& h(k) \int_{0}^{1} f_u^T(x(t(k)+\tau h(k),u_n) \lambda(t(k)+\tau h(k))  d\tau  \nonumber\\
&=& h(k) \sum_{i=1}^{s} b_i f_u^T \left(x \left(t(k)+c_i h(k)\right),u_n \right) \lambda \left( t(k)+c_i h(k) \right) + \nonumber \\
& & O(h^{\tilde{p}+1}(k)) . \nonumber
\end{eqnarray}
From (\ref{discrEq:errEstStateODEIntGlob}) and (\ref{discrEq:errEstAdjIntGlob}) we obtain
\begin{eqnarray}
&h(k)& \sum_{i = 1}^{s} b_i f_u^T(x_i(k+1),u_n)  \lambda_i(k) \\
&=& h(k) \sum_{i=1}^{s} b_i f_u^T \left(x \left(t(k)+c_i h(k)\right),u_n \right) \lambda \left( t(k)+c_i h(k) \right) \nonumber\\
&&+ O(h(k)h_m^{min\{ q+1, \bar{q}_d+1 \}}). \nonumber 
\end{eqnarray}
The local error of the discrete reduced gradients is therefore 
\begin{eqnarray}
&\int_{t(k)}^{t(k+1)}& f_u^T(x(t),u_n) \lambda(t)  dt - h(k) \sum_{i = 1}^{s} b_i f_u^T(x_i(k+1),u_n) \\
&=& O(h^{\tilde{p}+1}(k)) + O(h(k)h_m^{min\{ q+1, \bar{q}_d+1 \}}) =  \nonumber \\
&=& O(h(k)h_m^{min\{\tilde{p}, q+1, \bar{q}_d+1 \}}) = O(h(k)h_m^{\tilde{p}_d}) , \nonumber
\end{eqnarray} 
where $ \tilde{p}_d = min\{\tilde{p}, q+1, \bar{q}_d+1 \} $. For RADAU IIA scheme we obtain $ \tilde{p}_d= min\{5, 3+1, 2+1 \} = 3 $.
If the local error of the discrete reduced gradient is  $ O(h(k)h_m^{\tilde{p}_d}) $ then the global error can be estimated as  (\cite{hnw93}, Theorem 3.6 p. 162) 
\begin{equation}
\frac{d \phi(X(K))}{du_n} - \frac{d \phi(x(t_f))}{ du_n } =  O(h_m^{\tilde{p}_d}) 
\end{equation}
This result confirms that the discrete reduced gradients (\ref{discrEq:discrRedGrad002}) converge to the continuous reduced gradients (\ref{discrEq:redGradCont}) with an order at least $ \tilde{p}_d =3 $ (for RADAU IIA scheme). The discrete reduced gradients provide therefore an efficient and reliable method for approximation of reduced gradients.

\section{Conclusions}

The paper presents the novel computational approach to hybrid optimal control problems with sliding modes. Our computational method is based on a Runge--Kutta scheme for integrating system and adjoint equations. We show that in case of ODEs system equations, discrete adjoint equations correspond to the discretized continuous adjoint equations and on that basis we justify the convergence of the discrete adjoint equations to their continuous counterparts for step sizes approaching zero, we also give an estimate of the convergence order. 
In the second part of the paper we derive a similar result for index 2 DAEs describing the system motion in a sliding mode. In the second part we also present several numerical results of the application of the proposed computational procedure to optimal control problems with hybrid systems.


\begin{thebibliography}{99}
	
	\bibitem{ap1998}
	U. Ascher and L. Petzold,
	{\it Computer methods for ordinary differential equations and differential-algebraic equations },
	Society for Industrial and Applied Mathematics, Philadelphia, (1998).
	
	\bibitem{bh1975}
	A. Bryson and Y. Ho,
	{\it Applied Optimal Control },
	Hemisphere, New York, (1975).
	
	\bibitem{dl2009} 
	L. Dieci and L. Lopez, 
	"Sliding Motion in Filippov Differential Systems: Theoretical Results and a Computational Approach",
	{\it SIAM J. Numer. Anal.}, vol. 47(3), pp. 2023--2051, (2009). 
	
	\bibitem{h00} W. Hager, "Runge--Kutta methods in optimal control and the transformed adjoint equations", {\it Numer. Math.}, vol. 87, pp. 247--282, (2000).
	
	\bibitem{hlr89}
	E. Hairer, Ch. Lubich and M. Roche,
	{\it The Numerical Solution of Differential--Algebraic Equations by	Runge--Kutta Methods},
	Lecture Notes in Mathematics, vol. 1409, Springer--Verlag, Berlin, Heidelberg, (1989).
	
	\bibitem{hnw93}
	E. Hairer, S. Norsett and G. Wanner,
	{\it Solving Ordinary Differential Equations I},
	Springer--Verlag, Berlin, Heidelberg, (1993).
	
	\bibitem{hw96}
	E. Hairer and G. Wanner,
	{\it Solving Ordinary Differential Equations II},
	Springer--Verlag, Berlin, Heidelberg, (1996).
	
	\bibitem{ljs2000}
	J. Lygeros, K. H. Johansson, S. S. Sastry and M. Egerstedt,
	"On the existence of executions of hybrid automata",
	in {\it Proceed. of the 38th IEEE CDC, Phoenix, Arizona, 1999}, pp. 2249--2254, (1999).
	
	\bibitem{mhrl2015} I. Mynttinen, A. Hoffmann, E. Runge, P. Li, "Smoothing and regularization strategies for optimization of hybrid dynamic systems", {\it Optim Eng.}, Vol. 16, pp. 541--569, (2015).  
	
	\bibitem{pc2013} A. Pakniyat and P.E. Caines, "The Hybrid Minimum Principle in the Presence of Switching Costs ", {\it Proceedings of the 52nd IEEE Conference on Decision and Control, December 10-13, Florence, Italy}, pp. 3831--3836, (2013).
	
	\bibitem{pp02} O. Pironneau, E. Polak, "Consistent Approximations and Appropriate Functions and Gradients in Optimal Control", {\it SIAM J. Control and Optimization}, vol. 41, pp. 487--510, (2002).
	
	\bibitem{py99}
	R. Pytlak,
	{\it Numerical Methods for Optimal Control Problem With State Constraints},
	Lecture Notes in Mathematics, vol. 1707, Springer--Verlag, (1999), .
	
	
	\bibitem{py11} R. Pytlak, "Numerical procedure for optimal control of higher index DAEs", 
	{\it J. Discret. Dyn. Nat. Soc. (USA)}, vol. 29, pp. 1--24, (2011).
	
	\bibitem{ps19} R. Pytlak and D. Suski, "Algorithms for optimal control of hybrid systems with sliding motion",  arXiv:2101.04754, (2021).
	
	\bibitem{ps19b} R. Pytlak and D. Suski, "Numerical procedure for optimal control of hybrid systems with sliding modes, Part 2", submitted to arXiv, (2021).
	
	\bibitem{ps17} 
	R. Pytlak and D. Suski, 
	"On solving hybrid optimal control problems with higher index DAEs", {\it Optimization Methods \& Software} vol. 32, pp. 940--962, (2017).
	
	\bibitem{pz14} 
	R. Pytlak and T. Zawadzki, 
	"On solving optimal control problems with higher index DAEs", {\it Optimization Methods \& Software} vol. 29, pp. 1139--1162, (2014).
	
	\bibitem{ss2000}
	A. van der Schaft and H. Schumacher,
	{\it An Introduction to Hybrid Dynamical Systems},
	Springer-Verlag London, (2000).    
	
	\bibitem{s2004}
	M. S. Shaikh, "Optimal Control of Hybrid Systems: Theory and Algorithms", Ph.D. dissertation, McGill University, Montreal, May 2004. [Online]. Available: http://digitool.library.mcgill.ca/R/?func=dbin-jump-full\&object\_id=85095\&local\_base=GEN01-MCG02, (2004)
	
	\bibitem{sc2007} M. S. Shaikh, and P.E. Caines, "On the Hybrid Optimal Control Problem: Theory and Algorithms",
	{\it IEEE Transactions on AC, Vol. 52}, No. 9, pp. 1587--1603, (2007).
	
	\bibitem{sc2009} M. S. Shaikh, and P.E. Caines, "Correction to 'On the Hybrid Optimal Control Problem: Theory and Algorithms'", {\it IEEE Transactions on AC}, Vol. 54, No. 6, p. 1, (2009).
	
	\bibitem{sa2012} H. Schiller, M. Arnold, "Convergence of continuous approximations for discontinuous ODEs", 
	{\it Applied Numerical Mathematics}, Vol. 62, pp. 1503--1514, (2012).
	
	\bibitem{sp98} A. Schwartz and E. Polak "Consistent approximations for optimal control problems based on Runge--Kutta integration", {\it SIAM J. Control and Optimization}, vol. 34, pp. 1235--1269, (1998).
	
	\bibitem{sa2010} D. Stewart and M. Anitescu, 
	"Optimal control of systems with discontinuous differential equations", 
	{\it Num. Math.}, vol. 114, pp. 653-695, (2010).
	
	\bibitem{s1999}
	H. Sussmann,
	"A maximum principle for hybrid optimal control problems", 	{\it Proceed. of the 38th IEEE CDC, Phoenix, Arizona, 1999}, pp. 425--430, (1999).
	
	\bibitem{tc2010} F. Taringoo and P.E. Caines, "Gradient-Geodesic HMP Algorithms for the Optimization of Hybrid Systems Based on the Geometry of Switching Manifolds", {\it Proceedings of the 49th IEEE Conference on Decision and Control, December 15-17, Atlanta, GA, USA}, pp. 1534--1539, (2010).
	
	\bibitem{tc2011} F. Taringoo and P.E. Caines, "On the Extension of the Hybrid Minimum Principle to Riemannian Manifolds", {\it Proceedings of the 50th IEEE Conference on Decision and Control, December 12-15, Orlando, Fl, USA}, pp. 3301--3306, (2011).
	
	\bibitem{tc2012} F. Taringoo and P.E. Caines, "Newton-Geodesic HMP Algorithms for the Optimization of Hybrid Systems and the Geometric Properties of Hybrid Value Functions", {\it Proceedings of the 51st IEEE Conference on Decision and Control, December 10-13, Maui, Hawaii, USA}, pp. 4211--4216, (2012).
	
	\bibitem{tc2013} F. Taringoo and P.E. Caines, "On the Optimal Control of Hybrid Systems On Lie Groups and the Exponential Gradient HMP Algorithm", {\it Proceedings of the 52nd IEEE Conference on Decision and Control, December 10-13, Florence, Italy}, pp. 2653--2658, (2013).
	
	\bibitem{w1966}
	H. Witsenhausen, "A class of hybrid-state continuous-time dynamic systems", {\it IEEE Trans. on Autom. Control}, vol. 11, pp. 161--167, (1966).
	
\end{thebibliography}
\end{document}